\title{Twins and Vertex-Identification on Graphs}
\author{
        \textsc{(1)Sylvain Gravier}
            \qquad
      \textsc{(1)Simon Schmidt}
      \qquad  \textsc{(2)Souad Slimani}
        \thanks{Contact author}
       (1) Institut Fourier UGA-CNRS, SFR Maths \`a Modeler, \\100 rue des maths, 
       BP 74, 38402 St Martin d'H\`eres cedex, France.\\
  (2) Laboratoire LaROMaD. U.S.T.H.B Universit\'e, Facult\'e des Maths\\
BP 32 El Alia Bab Ezzouar 16 111 Alger}

\documentclass[12pt]{article}

\usepackage{amssymb}
\usepackage{amsmath}
\usepackage{amsthm}
\usepackage{color}
\usepackage{curves}
\usepackage{eepic}
\usepackage{epic}
\usepackage{graphics,color,xcolor}
\usepackage{pstricks}
\usepackage{pstricks-add,pst-node}

\newcommand{\Nb}{\mathbb{N}}

\newtheorem{theorem}{Theorem}[section]

\newtheorem{property}[theorem]{Property}
\usepackage{makeidx}              

\title{Twins and Vertex- Identification on Graphs}

   \author{Sylvain Gravier$^{1}$ \and Simon Schmidt$^{1}$  \and Souad  Slimani$^{1,2}$\\
  \small $^{1}$ Institut ~Fourier~ UGA-CNRS, SFR ~Maths~ \`a Modeler, ~100~ rue ~des ~maths,\\
\small BP 74, ~38402 ~St Martin~ d'H\`eres~ cedex, ~France \\
 \small  $^{2}$ Laboratoire~ LaROMaD. ~U.S.T.H.B ~Universit\'e,~ Facult\'e ~des ~Math\'ematiques~\\
\small BP 32 ~El Alia~ Bab ~Ezzouar ~16 111 Alger
}
\makeindex		    

\begin{document}

\maketitle




\begin{abstract}
\noindent Recently, several vertex identifying notions were introduced (identifying coloring, lid-coloring, ...), these notions were inspired by identifying codes. All of them, as well as original identifying code, are based on separating  two vertices according to some conditions on their closed neighborhood. Therefore, twins can not be identified. So most of known results focus on twin-free graph. Here, we show how twins can modify optimal value of vertex-identifying parameters for identifying coloring and locally identifying coloring.

\end{abstract}

\textbf{Keyword :}
 Identifying coloring, locally identifying coloring, twins, separating.
\section{Introduction}\label{sec1}
In this paper, we are interested on vertex-colorings allowing to distinguish the vertices of a graph. Given a graph $G$ and a coloring $c$, a pair of vertices $u$ and $v$ of $G$ is {\em identified} if and only if $c(N[u])\neq c(N[v])$ (where $N[u]$ denotes the closed neighborhood of $u$). Two vertices $u$ and $v$ with $N[u]=N[v]$ are called twins. Observe that two twins can not be identified.

\noindent During the discrete week of Institute Fourier at Grenoble in $2009$, Eric Duch\^{e}ne and Jullien Moncel presented the concept of {\em identifying coloring} of a graph as a vertex coloring such that any pair of vertices is identified.  Clearly, an identifying coloring of $G$ exists if and only if $G$ has no twins. 

\noindent Since, several authors \cite{ASS, lid1} introduced two notions of {\em locally identifying} where identified pair of  vertices is required only for adjacent vertices. Moreover, in order to incorporate graph with twins, the identifying condition concerns only pair of non-twin vertices.

\noindent Since only few results are known on identifying coloring (see \cite{these}), and in order to be uniform, we propose here to modify the definition of identifying coloring of a graph as a vertex coloring such that any pair of {\bf non-twin} vertices is identified.

\noindent In \cite{ASS}, a {\em relaxed locally identifying coloring}, {\em rlid-coloring} for short, of a graph $G=(V,E)$ is defined as a mapping $c:V\rightarrow \mathbb{N}$ such that any pair of {\bf adjacent} non-twin vertices is identified. Aline Parreau et al. \cite{lid1} introduced the notion of {\em locally identifying coloring}, {\em lid-coloring} for short, as a rlid-coloring $c$ which is proper that is $c(u)\neq c(v)$ for all pair of adjacent vertices $u, v$.

\noindent Given a graph $G$,  $\chi_{id}(G)$ (respectively $\chi_{lid}(G), \chi_{rlid}(G)$) denotes the smallest number of colors needed to have an identifying coloring (resp. lid-coloring, rlid-coloring) of $G$.

\noindent The $\chi_{lid}$ is the most studied of these parameters, see for instance \cite{lid1,lid2,lid3}. Nevertheless, except in \cite{ASS}, most of the results concern twin-free graphs. The aim of this paper is to show that twins may have significant influence on these parameters. 

\noindent In order to state our results, we will need additional definitions.
\noindent Let $\mathcal{R}$ be the equivalence relation defined as follows: for all vertices $u,  v \in V(G)$, we have $u \mathcal{R}v $ iff $N[u]=N[v]$.
\noindent Denote by $G \backslash\mathcal{R}$, the maximal twin-free subgraph of $G$ (that is the quotient of $G$ by relation $\cal R$). The number of  equivalence-classes having at least two vertices in $G$ is denoted by $t(G)$. Let $T(G)$ be the cardinality of a largest equivalence-class.

\noindent In \cite{ASS}, the authors proved the following theorem~:

\begin{theorem}\label{thm:SS}
Let $G$ be a graph. Then we have $$\chi_{rlid}(G\backslash{\cal R})-t(G)\leq \chi_{rlid}(G)\leq \chi_{rlid}(G\backslash{\cal R}).$$
\end{theorem}

\noindent Moreover, the authors in \cite{ASS} exhibit graphs for which the bounds are tight. In this paper, we give analogous results for identifying colorings and lid-colorings.

\begin{theorem}\label{thm:id}
Let $G$ be a graph. Then we have $$\chi_{id}(G\backslash{\cal R})-t(G)\leq \chi_{id}(G)\leq \chi_{id}(G\backslash{\cal R}).$$
\end{theorem}

\begin{theorem}\label{thm:lid}
Let $G$ be a graph. Then we have $$\chi_{lid}(G\backslash{\cal R})-t(G)\leq \chi_{lid}(G)\leq \chi_{lid}(G\backslash{\cal R})+(T(G)-1).t(G).$$
\end{theorem}

\noindent Proofs of Theorems \ref{thm:id} and \ref{thm:lid} are given in Section \ref{sec:proofs}. In Section \ref{sec:extremal}, we exhibit graphs for which the bounds in Theorems \ref{thm:id} and \ref{thm:lid} are tight.

\section{Proofs of bounds}\label{sec:proofs}

\noindent {\bf Proof of Theorem \ref{thm:id}.}
\noindent We present here a proof similar to the proof of Theorem \ref{thm:SS} given in \cite{ASS}.

\noindent Consider an identifying coloring $c$ of $G\backslash{\cal R}$. Let define a coloring $c'$ of $G$ as follows~: for each vertex $x$ in $G\backslash{\cal R}$ and its twin (if there exists) $y$, set $c'(x)=c'(y)=c(x)$.
\noindent Since in $G$, we are not interested to distinguish the twins then $c$ defines an identifying coloring of $G$.

\noindent Now, let $c$ be an identifying coloring of $G$ using colors $\{1,\ldots, \chi_{id}(G)\}$. Consider the coloring $c'$ defined as follows : 
$c'(u)=c(u)$ if the vertex $u$ has no twin in $G$ and 
color the other $t(G)$ vertices of $G\backslash{\cal R}$ with different colors $\chi_{id}(G)+1$ until $\chi_{id}(G)+t(G)$. 
\noindent This coloring gives an identifying coloring of $G\backslash{\cal R}$. \qed

\vspace{0.5cm}
\noindent {\bf Proof of Theorem \ref{thm:lid}.}\\
\noindent Now, consider a $lid$-coloring $c$ of $G\backslash{\cal R}$ using colors $\{1,\ldots, \chi_{lid}(G\backslash{\cal R})\}$. 
\noindent Let $c'$ be a coloring obtained from $c$ as follow: $c'(u)=c(u)$ for all $u$ in $G\backslash{\cal R}$. By definition, there are at most $(T(G)-1).t(G)$ vertices in $G$ which are not in $G\backslash{\cal R}$. For each of them assign a distinct color from $\{\chi_{lid}(G\backslash{\cal R})+1,\ldots, \chi_{lid}(G\backslash{\cal R})+(T(G)-1)t(G)\}$. This coloring gives an lid-coloring of $G$.

\noindent Now, similarly as proof of Theorem \ref{thm:id}, let $c$ be a $lid$-coloring of $G$ using colors $\{1,\ldots, \chi_{lid}(G)\}$. Consider the coloring $c'$ defined as follows : 
$c'(u)=c(u)$ if the vertex $u$ has no twin in $G$ and 
color the other $t(G)$ vertices of $G\backslash{\cal R}$ with different colors $\chi_{lid}(G)+1$ until $\chi_{lid}(G)+t(G)$. 
\noindent This coloring gives an lid-coloring of $G\backslash{\cal R}$. \qed

 \section{Extremal graphs}\label{sec:extremal}
 
\noindent  First define the split graph $H_p= (S_p\cup K_p, E)$ for a given integer $p$ where $S_p=\{s_1, \ldots, s_p\}$ (respectively $K_p=\{k_0,\ldots, k_p\}$) induces a stable (resp. clique). The others edges of $H_p$ are $s_ik_i$ for all $i=1, \ldots, p$.
 
\begin{property}\label{pr:Hp}
Let $p\geq 1$ be an integer. We have that $$\chi_{id}(H_p)= p+2\hbox{ and } \chi_{lid}(H_p)= 2p+1.$$
\end{property}

\begin{proof}
\noindent The coloring $c$ defined by $c(s_i)=i, c(k_i)=p+1$ for all $i=1,\ldots p$ and $c(k_0)=p+2$, is an identifying coloring of $H_p$.\\
\noindent Let prove now that $\chi_{id}(H_p)\geq p+2$. Let $c$ be an identifying coloring of $H_p$. First observe that $c(s_i)\neq c(s_j)$ for all $i\neq j$. Indeed, otherwise $c(N[k_i])=c(N[k_j])$ which leads a contradiction. Second, suppose that $c(s_i)=c(k_j)$ for some integers $i, j$ ($i$ can be equal to $j$). Then $c(N[k_0])=c(N[k_i])$, a contradiction. To conclude, check that if $\vert c(K)\vert =1$ then $c(N[s_i])=c(N[k_i])$ for all $i$, a contradiction which completes the proof of $\chi_{id}(H_p)\geq p+2$.

\noindent Any coloring using $2p+1$ distinct colors is a lid-coloring of $H_p$.\\
\noindent Let prove now that $\chi_{lid}(H_p)\geq 2p+1$. Let $c$ be an lid-coloring of $H_p$. As previously,  we have $c(s_i)\neq c(s_j)$ for all $i\neq j$ else $c(N[k_i])=c(N[k_j])$ and $k_i$ and $k_j$ are adjacent. Again, $c(s_i)\neq c(k_j)$ for all pair $i, j$, otherwise $c(N[k_0])=c(N[k_i])$ for some $i\neq 0$, a contradiction. To conclude, since $K$ is a clique, then $\vert c(K)\vert = p+1$. 
\end{proof}
 
 \vspace{0.5cm}
\noindent Consider the first extension $H^{ext}_{2^a}= (S_{a.2^a}\cup K_{2^a}, E)$ for some integer $a\geq 1$ where $K_{2^a}$ induces a clique. One may define the vertices of $K_{2^a}= \{k_{\cal E} \vert {\cal E}\subseteq \{1,\ldots, a\}\}$. Now define $S_{a.2^{a-1}}=\{s_{{\cal E}_i} \vert {\cal E}\subseteq \{1,\ldots, a\} \hbox{ and } i\in {\cal E}\}$. Observe that $\vert K_{2^a}\vert = 2^a$ and $\vert S_{a.2^{a-1}}\vert=a.2^{a-1}$. The others edges of $H^{ext}_{2^a}$ are $s_{{\cal E}_i} k_{\cal E}$ for all $i\in {\cal E}$ and $s_{{\cal E}_i}s_{{\cal E}_j}$ for all $i, j\in {\cal E}$.\\
Remark that $H_{2^a}^{ext}\backslash{\cal R}=H_{2^a-1}$ with $t(H_{2^a}^{ext})=2^a-1-a$ and $T(H_{2^a}^{ext})=a$.\\

\begin{property}\label{pr:Hpext}
Let $a\geq 1$ be an integer. We have that $$\chi_{id}(H_{2^a}^{ext})= a+2\hbox{ and } \chi_{lid}(H_{2^a}^{ext})= a+2^a.$$
\end{property}

\begin{proof}
\noindent The coloring $c$ defined by $c(s_{{\cal E}_i})=i, c(k_{\cal E})=a+1$ for all ${\cal E}\subseteq \{1,\ldots, a\}$ and for all $i\in {\cal E}$ and $c(k_0)=a+2$, is an identifying coloring of $H_{2^a}^{ext}$.\\
\noindent By Theorem \ref{thm:id}, we have $\chi_{id}(H_{2^a}^{ext})\geq \chi_{id}(H_{2^a}^{ext}\backslash{\cal R})-t(H_{2^a}^{ext}) = 2^a -1 + 2 - (2^a -1 - a) = a+2$.

\noindent For the lid-coloring the proof is similar except that we need distinct colors for each vertex in the clique $K_{2^a}$.  
\end{proof}
 \vspace{0.5cm}

 \noindent The two previous propositions show that lower bound of Theorems \ref{thm:id} and \ref{thm:lid} are tight for graph $H_{2^a}^{ext}$. Upper bound of Theorem \ref{thm:id} is reached for any twin-free graph. 
 
 \noindent Given integers $p\geq  t\geq 1$ and $T\geq 1$, consider the graph $H_p^{(T, t)}$ obtained from $H_p$ by adding $T-1$ twins to all vertices  $k_i$ for all $i=1, ...., t$. \\
 \noindent Remark that $H_{p}^{(T, t)}\backslash{\cal R}=H_{p}$ with $t(H_{p}^{(T, t)})=t$ and $T(H_{p}^{(T, t)})=T$.\\

  \begin{property}\label{pr:HTt}
Let $a\geq 1$ be an integer. We have that $$\chi_{lid}(H_{p}^{(T, t)})= 2p+1+(T-1).t.$$
\end{property}

\begin{proof}
\noindent Any coloring using $2p+1+(T-1).t$ distinct colors is a lid-coloring of $H_p^{(T, t)}$.\\
\noindent The proof of $\chi_{lid}(H_{p}^{(T, t)})\geq 2p+1+(T-1).t$ follows the one of Proposition \ref{pr:Hp}.  
\end{proof}
 \vspace{0.5cm}
\noindent For all triple of integers $p\geq  t\geq 1$ and $T\geq 1$, Proposition \ref{pr:HTt} show that upper bound of Theorem \ref{thm:lid} is reached.
 
 \section*{Concluding Remarks}
\noindent The main motivation of the present paper is to point out that twins may play a crucial role in identifying coloring problems using closed neighborhood. Probably it should be to difficult to re-consider all known results on twin-free graphs. \noindent But there are some special classes of graphs (e.g. split graphs) for which this work remains attractive.\\
\noindent Instead of coloring, one may ask what happens for identifying codes~? An identifying code \cite{karpovsky} is a subset of vertices $C$, such that for any pair of distinct vertices $u, v$, $N[u]\cap C\neq N[v]\cap C$.  This is the classical definition, and clearly, a graph having twins does not admit an identifying code. Consider, now the new definition where the condition  $N[u]\cap C\neq N[v]\cap C$ has to be verified only for non-twin pair of vertices. \\
It is not too difficult to see that the size of a minimum identifying codes in a graph $G$ with new definition is equal to the size of a minimum identifying codes in $G\backslash{\cal R}$. Therefore, it is not restrictive for identifying codes to consider only twin-free graphs.\\
\noindent Now, for coloring versions of identifying problems one may consider a {\em weighted} version. Given a graph $G$ and a weight function $w:V\rightarrow \Nb$, a mapping $c:V\rightarrow {\cal P}(\Nb)$ is an {\em weighted-identifying coloring} of $G$  if and only if $\vert c(u)\vert \leq w(u)$ for all vertices $u$ and all distinct pairs of non-twin vertices are identified.\\
Observe that an optimal value for a pair $(w, G)$ is the same than the optimal value for the pair $(w', G\backslash{\cal R})$ where $w'(u)=w(u)+T(u)-1$ where $T(u)$ is the number of twins of $u$.  So for this new definition one may focus only on twin-free graphs.

\end{document}